# Topología Algebraica para el Análisis de Datos

Karla Saraí Jiménez-Martínez , Daniel Trejo-Medina
DIVULGACIÓN CIENTÍFICA DSA SOLUCIONES®
Ciudad de México

___________________________________________________________________

**Resumen - Abstract**

En esta investigación se aborda una novedosa herramienta para el análisis de datos, conocida como Análisis Topológico de Datos (TDA por sus siglas en inglés). Subyace de un área de las Matemáticas conocida como Álgebra Combinatoria o más recientemente Topología Algebraica, la cual a través de hacer un fuerte uso de Computación, Estadística, Probabilidad y Topología, entre otros conceptos, extrae características matemáticas de un conjunto de datos que nos permiten asociar, crear e inferir información general y de calidad sobre estos.

**Palabras clave:** Topología Algebraica, Topología Combinatoria, TDA, Análisis topológico de datos.

1. **Introducción**

Con el advenimiento y uso de nuevas tecnologías de información, entre ellos dispositivos móviles, documentos o datos creados en redes sociales y su comentarios adicionales, han generado una cantidad enorme de datos.

Estos datos sin un propio uso, carecen de valor, sin embargo, imagine que es un banco y desea disminuir el riesgo de atraer un nuevo cliente, ya sea por evitar un posible lavado de dinero, o una aseguradora que busca mitigar el riesgo utilizando un análisis de riesgo no estructurado y de fuentes informales.

De inicio usted puede considerar que este conjunto de datos no es factible de analizar y generar o inferir información general y sobre todo con calidad para poder tomar una decisión.

Tradicionalmente los métodos bayesianos son los mas comunes y conocidos en la jerga de la computación por la popularización de diversas herramientas que lo tienen como base; con el manejo de contenido (por ejemplo videos de YouTube, imágenes de Instagram, noticias de periódicos, reportes de buholegal.com), datos no estructurados y las tradicionales bases de datos (no SQL y SQL) hemos utilizado otras aproximaciones para el análisis de datos, donde herramientas como Raytheon, TIBCO Spotfire, Attivio, GreenplumDB, R y otras permiten su aplicación y uso real.

Si es usted un banco que necesita medir riesgos operativos, o prevenir fraude interno, es un vendedor al detalle (*retailer*) que busca como optimizar el valor agregado a sus clientes y mitigar las mermas sociales y de baja de marca, o un órgano de seguridad publica que busca analizar múltiples factores en tiempo corto, una telefónica esperando ajustar el nivel de tasa de abandono en tiempo real desde múltiples dimensiones, este artículo le puede interesar.

En un ámbito mas directo en recaudación de impuestos es el cruce de causantes menores, mayores, facturas y transferencias bancarias, es decir la recaudación de impuestos donde se puede aplicar de forma mas directa.

El tener múltiples fuentes de datos, permite mediante el uso de espacios topológicos euclidianos, el generar una nube de ellos que forman un conocimiento generalizado sobre ellos, lo cual da como resultado una vinculación y escenario de posibles pronósticos de datos que pueden ayudar a resolver diversas preguntas de negocio.

La visualización de datos de demasiadas dimensiones (para fines prácticos mas de 25) genera una dificultad de visualización directa, sin embargo, en la aplicación real, no exceden de 10 dimensiones las que se analizan.

En las últimas décadas, la cantidad de datos que se crea y comparte en Internet ha tenido un crecimiento asombroso, tan sólo las transacciones financieras, GPS y redes sociales, generan 2.5 quintillones de exabytes cada día alrededor del mundo y se espera que la cantidad siga aumentando. Para analizar todos estos datos que los individuos, empresas públicas y privadas o cualquier ente genera en grandes cantidades, han surgido disciplinas que buscan lograr un eficiente manejo, tal es el caso del Análisis Topológico de Datos, área de aproximadamente 15 años de antigüedad, proveniente de la Topología Algebráica, subdisciplina de las Matemáticas.

A continuación presentamos una explicación que consideramos sencilla de los métodos que ocupa el DSAengine® para uno de sus modelos analíticos aplicados.

La manera en que opera este tipo de análisis, es partiendo de una base o conjunto de datos, del cual se desea extraer características que nos permitan generar conocimiento generalizado sobre ellos. Los datos se imaginan como una nube, formada por puntos, en los que cada uno representa un dato. A esta nube se le asocia un espacio topológico euclidiano, conocido con el nombre de Complejo Simplicial, este se vuelve en el interés de estudio para extraer sus características topológicas, tales como sus componentes conexos, agujeros, su estructura como gráfica, ente otras. Asociar una estructura topológica a los datos permitirá pensar en si esta se parece a algún objeto matemático del que ya se conocen sus propiedades topológicas y así poder trasladarlas a estos para su estudio.

Cada punto que representa un dato en el espacio, se toma como centro para trazar un círculo de diámetro d fijo, aquí resulta claro que como los datos están dispersos de manera aleatoria, hay puntos muy cercanos y lejanos unos de otros, encontrando círculos que se intersectan y otros no. Así podrán intersectarse, ninguno, dos, tres, cuatro…n círculos. *Figura 1*.

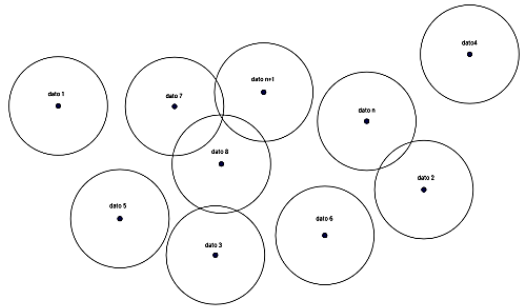

**Figura 1**. Comportamiento de los círculos trazados a partir de cada punto que representanta un dato.

Cuando se intersectan dos círculos, se traza un vértice que une a los dos centros de estos, si se intersectan tres, se trazará un triángulo, para cuatro un tetraedro, y así sucesivamente. Al realizar este procedimiento con cada uno de los datos, el resultado final será una gráfica formada con todos los vértices, la cual incluirá hoyos o agujeros entre los círculos como en la *figura 2*. La gráfica se vuelve en el espacio topológico de estudio, del cual se desea conocer y analizar sus propiedades, principalmente los agujeros presentes en ella, ya que la presencia de estos en el tiempo, permiten conocer información sobre la base de datos.

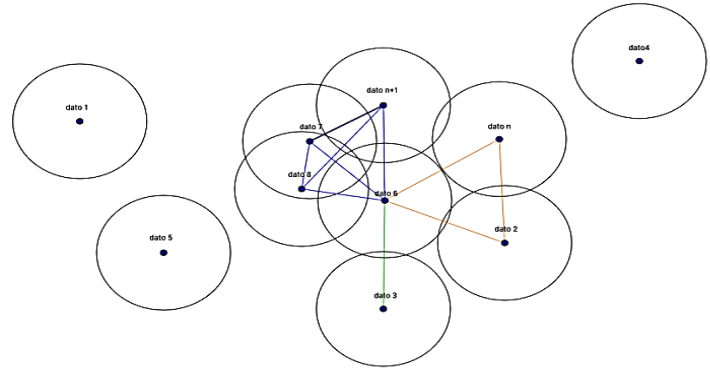

**Figura 2**. Gráfica formada por los vértices de los círculos que se intersectan.

Al estudiar los agujeros presentes en la gráfica, es necesario establecer un diámetro $d$ fijo adecuado con el que se trazarán los círculos, ya que de este dependerá si los agujeros son grandes o chicos. Por ejemplo, si se supone un diámetro muy pequeño, los círculos alrededor de cada punto no se intersectarán y por el contrario, si se elige uno muy grande, todos los círculos estarán encimados o demasiado pegados ya que se intersectarán una gran cantidad, lo cual en una base de datos podría ser equivalente a pensar que se tienen datos atípicos, que se desvían mucho de los demás o la presencia de ruido. *Figura 3*.

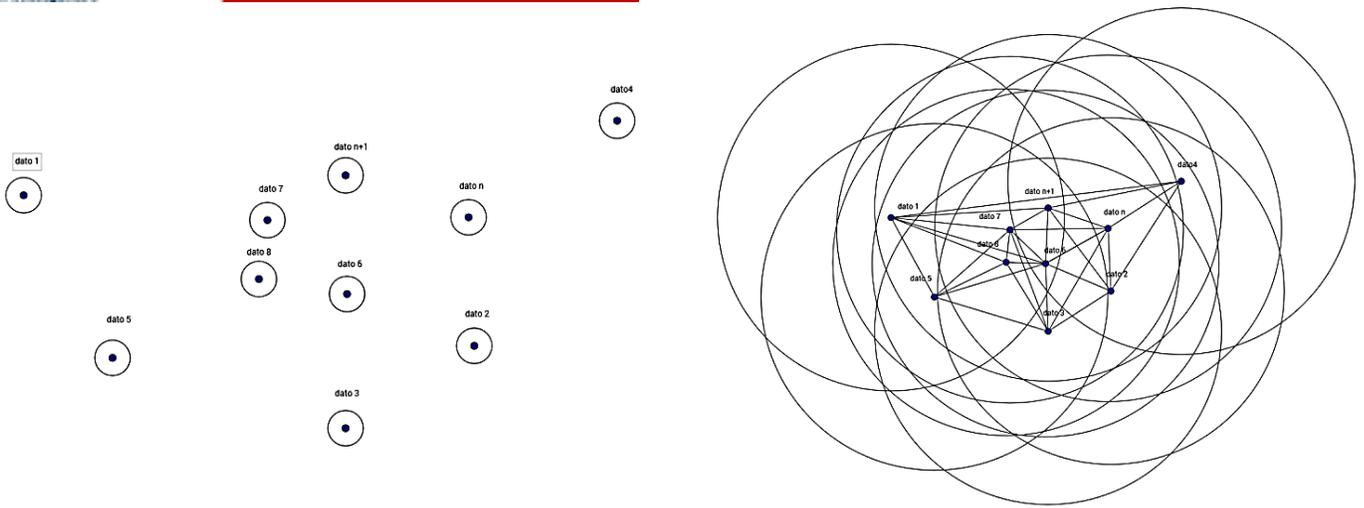

**Figura 3.** Comparación de elegir un diámetro d muy pequeño, en el que ningún círculo se intersecta (izq) y elegir un diámetro d muy grande, en el que todos los círculos se intersectan (der)

Ahora, suponga que se forma un agujero entre 4 círculos de diámetro $d_1$ provenientes de 4 puntos, y se busca desaparecer ese hoyo, necesariamente para lograrlo se tiene que modificar el diámetro aumentándolo, es decir se tendría un nuevo $d_2$. Los vértices que unen los puntos permanecerían, pero el hoyo quedaría oculto, como en la *figura 4*.

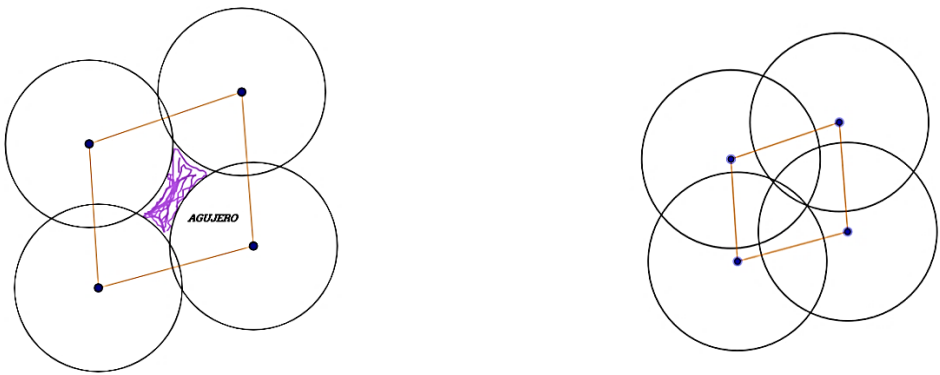

**Figura 4.** Presencia y desaparición del agujero, conforme se modifica el diámetro $d$. Los vértices permanecen.

Se vuelve entonces de especial interés qué tan persistente es o no un hoyo en el tiempo y su impacto en la nube de datos conforme se va modificando el diámetro de los círculos, para ello es que ser recurre al uso de la homología persistente. La persistencia de un hoyo en la gráfica se puede asociar a un segmento formado por la pareja $(d_1, d_2)$ partiendo del diámetro inicial $d_1$ para los círculos en la que el hoyo está presente, y posteriormente despareciendo de la gráfica, cuando $d_1$ se modifica y se convierte en un nuevo diámetro $d_2$. Cada agujero presente tendrá asociado su propio segmento $(d_1, d_2)$, formando entre todos estos un conjunto que se conoce con el nombre de código de barras, es decir que una nube de datos diámetros, tendrá asociado una identificación, un único código de barras.

Los hoyos que persisten mayoritariamente a través de la modificación del diámetro, y por tanto tienen un segmento o barra larga, son los de especial interés ya que aportan información y significado de la estructura de la base de datos. Por último, todo el conjunto de barras se transforma en puntos, para formar un diagrama de persistencia, el cual es una herramienta gráfica que nos permite resumir toda la información de los datos de estudio.

En la siguiente sección se describen los conceptos en los cuáles se apoya la Topología Algebraica para el análisis de datos, además de dar estos en el orden en el que van operando sobre la nube de datos para extraer las características que se desean. Se omiten demostraciones ya que no es el objetivo de un artículo de divulgación, sin embargo en la parte última, se sugiere la bibliografía al respecto.

2. **Conceptos**

En 1894, Henri Poncairé, a través de publicar varios artículos, comienza a sentar las bases de la Topología Algebraica. Creó lo que hoy se conoce como homología simplicial; describiendo los conceptos de triangulaciones, complejo dual, números de Betti, entre otros. Cada vez se han descubierto dentro y fuera de las matemáticas diversas aplicaciones de estos conceptos y actualmente, se han producido métodos para calcular información homológica a partir de una nube de datos X. Para lo cual es necesario el conocimiento de diferentes conceptos puramente matemáticos. En este apartado se darán definiciones básicas de topología algebraica, la cual estudia los invariantes respecto a funciones continuas.

**Topología**

**Def.2.1.** Sea X un conjunto no vacío y $\tau$ una familia de subconjuntos de X. Diremos que $\tau$ es una topología en X si satisface las siguientes condiciones:

- $\emptyset, X \in \tau$
- Si $A_1, A_2, \ldots A_n \in \tau \rightarrow \bigcap_{i=1}^{n} Ai \in \tau$
- $\{A_i\}_{i \in I}$

A la pareja (X, $\tau$) se le conoce como un espacio topológico. Se le llama puntos a los elementos de X y a los elementos de $\tau$, conjuntos abiertos. Si (X, $\tau$) es un espacio topológico, diremos que A⊆ X es un conjunto cerrado si X/A es un conjunto abierto. Las dos topologías más elementales son la topología discreta, la cual consta de todos los subconjuntos de un conjunto X y la topología indiscreta, que consta del vacío y del mismo conjunto X (i.e, $\tau$ = {$\emptyset$, X}).Esto quiere decir que cualquier conjunto X admite al menos dos topologías.

**Def.2.2** Sean $(X, \tau)$ y $(Y, \tau')$ dos espacios topológicos y F una aplicación entre ellos, diremos que F es continua si para cada conjunto abierto $A \in \tau'$, $F^{-1}(A)$ es un conjunto abierto en $\tau$, donde

$$F^{-1}(A) = \{x \in X | F(x) \in A\}$$

Para hablar de dos formas equivalentes en topología, la noción de equivalencia es la de homeomorfismo.

**Def. 2.3**. Se dice que dos espacios topológicos $(X, T)$ y $(X', T')$ son homeomorfos (o equivalentes) si existe una función biyectiva continua y con inversa continua. A una función que cumple estas condiciones se le conoce como homeomorfismo.

Si el interés es mostrar que dos espacios topológicos son homeomorfos, se deben encontrar propiedades en ellos que se preserven bajo homeomorfismos, así si uno de los espacios posee una de estas propiedades y el otro no, se concluye que no son homeomorfos.

**Def.2.4**. Se dice que una aplicación $\varphi : A \to B$ es un morfismo si para cada par de elementos x, y $\in A$ se tiene que $\varphi(x * y) = \varphi(x) *' \varphi(y)$, donde $*$ y $*'$ son las operaciones algebraicas en A y B respectivamente.

**Def.2.5.** Una categoría es una clase de objetos.

Por ejemplo, en la clase de los grupos abelianos, los objetos son los grupos abelianos. Dos categorías se relacionan mediante aplicaciones llamadas morfismos, en el ejemplo de los grupos, dichos morfismos son los homomorfismos de grupos. Si construimos ahora un universo de categorías este también será una categoría, y se pueden relacionar sus objetos (que son las categorías) mediante funtores los cuales se ven como una generalización del concepto de función para categorías. Los funtores asocian a cada objeto de una categoría, un objeto de la otra, y a cada aplicación de la primera una aplicación de la segunda. Para el Análisis Topológico de Datos, se usa el funtor de homología de complejos simpliciales, que a continuación se describe.

**Complejo simplicial**

Los complejos simpliciales son estructuras combinatorias que permiten la intervención del Álgebra en la Topología gracias a homología simplicial y los invariantes asociados a ella. La noción de complejo simplicial se desarrolló gradualmente a partir de los estudios de polígonos y poliedros tridimensionales que se remontan a los orígenes de las Matemáticas. (Elementos de la Homología Clásica Notas Teóricas, 2012)

**Def.2.6**. El complejo simplicial abstracto y finito, es una familia K no vacía de subconjuntos de un conjunto de vértices $V = \{v_i\}^m$ tal que se cumplen:

- Si $\alpha \in K$ y $\beta \subseteq \alpha$, entonces $\beta \in K$;
- $V \subseteq K$ (Para simplificar la notación, se identifica a cada $v \in V$ con $\{v\} \in K$)

A los elementos de K se les denomina caras, y su dimensión se define como uno menos que su cardinalidad. Las caras de dimensión cero se les denomina vértices y a las de dimensión uno aristas. Un mapeo simplicial entre complejos simpliciales es una función que respeta sus contenidos estructurales al mapear caras de una estructura a caras de la otra, es decir que estos conceptos presentan estructuras combinatorias que capturan las propiedades topológicas de una gran variedad de estructuras geométricas.

Dado un complejo simplicial abstracto K, se produce un espacio topológico simplicial, pues se considera la realización geométrica o poliedro asociada, que se denota como [K], la cual es construida cuando se consideran las caras de K como generalizaciones de triángulos y tetraedros en espacios euclidianos de alta dimensión.

Para analizar un complejo simplicial K, se construyen estructuras algebraicas y se calculan sus invariantes topológicos, que son propiedades de K y que no cambian bajo homeomorfismos, como se definió en el apartado anterior. Es decir que se calculan los invariantes topológicos de K al "traducir" su estructura combinatoria al álgebra, considerando el siguiente procedimiento.

1) Se construye un módulo de k cadenas $C_k$, que será resultado de todas las combinaciones lineales de caras k-dimensionales de K, con coeficientes dentro de un anillo conmutativo.

2) Se consideran los operadores frontera (o simplemente fronteras) $\partial_k: C_k \to C_{k-1}$ que son los morfismos que mandan una cara del conjunto $\sigma = [p_0, p_1, \ldots p_k] \in C_k$ en

$$\partial_k(\sigma)\sigma = \sum_{i=0}^{k} (-1)^i [p_0, \ldots p_{i-1}, p_{i+1}, \ldots p_k]$$

3) Se construye el grupo de homología de nivel k, definido por los módulos de cociente $H_k(K) \coloneqq \ker(\partial_k)/im(\partial_{k+1})$. Así entonces, se define el número de huecos k-dimensionales o k-ésimo número de Betti de K como $\beta_k = \text{rango}(H_k)$. Es decir que si tomamos de ejemplo la esfera, se tienen cero huecos 1-dimensionales y un solo hueco 2-dimensionales.

**Def.2.7.** Los operadores de frontera $\partial_k$ denotan un complejo de cadenas denotado por $C_* = C_*(K)$ representado por

$$\cdots \to C_{k+1} \xrightarrow{\partial_{k+1}} C_k \xrightarrow{\partial_k} C_{k-1} \to \cdots$$

**Def. 2.8.** [4] Dado un complejo de cadenas $C_*$ de módulos sobre un anillo conmutativo y unitario R, se definen los módulos de k-ciclos y k-fronteras como $Z_k = \ker \partial_k$ y $B_k = im \partial_{k+1}$, respectivamente. Como se tienen sub módulos anidados $B_k \subseteq Z_k \subseteq C_k$ el R-módulo de k-homología $H_k = (C_*) = Z_k/B_k$ está bien definido.

Durante todo el análisis Topológico de Datos, uno de los principales intereses es cuantificar los agujeros de K, cabe la posibilidad de encontrarse con dos ciclos que representen el mismo agujero, lo cual indicaría que no hay ningún espacio entre estos dos cíclos, ya que de existir dicho espacio saldría de entre ellos.

**Filtraciones**

**Def. 2.10**. Sea K un complejo simplicial de dimensión finita n, una filtración de K es una colección $K_0, K_1 \ldots K_m$ de complejos, tales que:

- $K_1 \subset K_2 \subset \cdots \subset K_m = K$

- $K_i$ es un subcomplejo de $K_{i+1}$ para i= 0,1…,m-1

El objetivo más importante es que conforme va formándose el complejo, se encuentren filtraciones con características que puedan aportar información importante para el análisis de la nube de datos Figura, con esto es natural pensar que no existirá una única filtración para un complejo K.

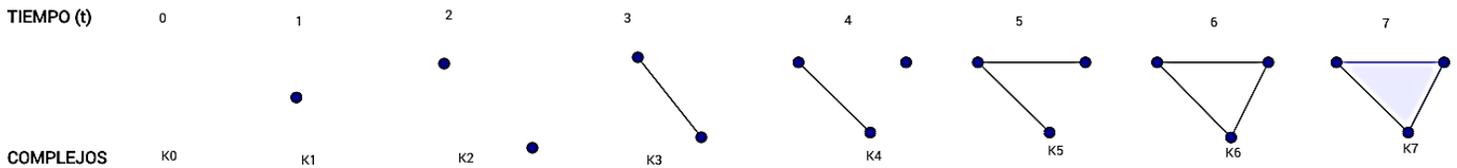

**Figura 5**. Filtración ordenada de un complejo simplicial formado por un triángulo y sus caras

**Homología persistente**

Ahora, situándose en el conjunto de datos denotado por $X = \{X_i\}_{i=1}^{m} \subset \mathbb{R}^n$ del que se desea obtener información, como se ha mencionado, los invariantes topológicos son características importantes de objetos geométricos y sus propiedades interesan ya que son los indicadores para entender al conjunto de datos. Al calcular estos invariantes, un problema a tratar es la inestabilidad que presenta la información homológica, ya que pequeñas variaciones al construir estructuras topológicas sobre X, podrán producir grandes cambios homológicos, para ello se usa la homología persistente, que permite calcular a través de herramientas computacionales, los invariantes topológicos en estructuras finitas. La homología es un modo de contar agujeros de cualquier dimensión en un espacio topológico, y en ese sentido, es una medida de la complejidad de dicho espacio. Su objetivo principal es entender cómo se relacionan distintos grupos de homologías por medio de inclusiones. (González, J. y Guillemard, M)

Intentar construir un complejo simplicial K, a partir de una nube de datos X puede ser un problema difícil, pues dependerá del número de datos con los que se cuenten en dicha base de datos. Se puede considerar la homología $X_\varepsilon = \cup_{i=1}^m B(x_i, \varepsilon)$ donde cada B representa una bola de radio $\epsilon$, que se traza tomando como punto central cada dato $x_i$, en este caso se buscaría un valor óptimo para el trazo de las bolas en el que la homología de $X_{\varepsilon 0}$ corresponde a la homología de una subvariedad M de $\mathbb{R}^n$ con lo que obtendríamos valores homológicos distintos por cada pequeña variación que se haga para $\varepsilon_0$ y con ello gran inestabilidad. Para resolver este problema, la homología persistente no considera información topológica para exclusivamente un valor fijo de $\varepsilon$ sino para todo $\varepsilon > 0$.

Gracias al concepto de filtración es que se puede asegurar que habrá un número finito de complejos simpliciales no homeomorfos $K_1 \subset K_2 \subset \cdots \subset K_m$ construidos a partir de $\{X_\varepsilon, \varepsilon > 0\}$ de los que se desea estudiar su información homológica, para lo cual existen varias estructuras simpliciales, en particular la construcción computacional descrito a continuación.

**Def. 2.11**. Complejo de Vietoris-Rips $R_\varepsilon(X)$,

- X se ve como un conjunto de vértices

- El conjunto de vértices $\sigma = \{x_0, \ldots x_k\}$ determina un k-simplejo de $R_\varepsilon(X)$ si $d(x_i, x_j) \leq \varepsilon$ para todo $x_i, x_j \in \sigma$

- Dado un valor $\varepsilon_i$, el complejo de Vietoris-Rips $R_\varepsilon(X)$, determina un elemento de la filtración $K_1 \subset K_2 \subset \cdots \subset K_m$ con $K_i = R_{\varepsilon i}(X)$

Es decir que un número finito de valores $\{\varepsilon_i\}_{i=1}^m$ describirán las características homológicas de X, formando con cada valor un complejo de Vietoris-Rips $K_i$ que resultará en una colección que representa las propiedades topológicas de la familia $\{X_\varepsilon, \varepsilon > 0\}$. Por lo tanto, el análisis topológico de una nube de puntos X y el crecimiento de su K complejo asociado se reduce al análisis de una filtración $K_1 \subset K_2 \subset \cdots \subset K_m$ y las clases de homología, principal objeto de estudio de la homología persistente.

**Persistencia**

Al querer estudiar la persistencia del complejo K asociado a un conjunto de datos X a través del tiempo, se busca poder tener información de generación y persistencia de las clases, respondiendo a las preguntas ¿Cuánto viven? ¿En qué momento es que nacen y mueren? Y, posteriormente resumir esta información en un gráfico llamado diagrama de persistencia, siendo necesarias las siguientes definiciones.

**Def. 2.12.** Complejo persistente. Familia de complejos de cadenas $\{C_*^i\}_{i \geq 0}$ y sus morfismos:

$$C_*^0 \xrightarrow{f^0} C_*^1 \xrightarrow{f^1} C_*^2 \xrightarrow{f^2} \cdots \xrightarrow{f^{i-1}} C_*^i \xrightarrow{f^i} C_*^{i+1} \xrightarrow{f^{i+1}} \cdots$$

Dada la filtración de K, se pueden considerar las funciones $f_i$ o inclusiones entre cada complejo simplicial de la sucesión $K_1 \subset K_2 \subset \cdots \subset K_m = K$

**Def.2.13.** Modulo persistente. Familia de R-módulos $M^i$ y homomorfismos $\varphi^i: M^i \dashrightarrow M^{i+1}$, el cual será de tipo finito si cada $M^i$ es finitamente generado y los mapeos $\varphi^i$ son isomorfismos para $i$ suficientemente grande.

**Def 2.14.** Sean $Z_k^i$ los módulos de k-ciclos y $B_k^i$ los módulos de k fronteras en la cadena $C_i$, los módulos p- persistentes de homología están definidos como:

$$H_k^{i,p} = Z_k^i \Big/ (B_k^{i+p} \cap Z_{k^i})$$

Donde el rango de $H_k^{i,p}$ es el k-ésimo número p-persistente de Betti de $C^i$ denotado por $\beta_k^{i,p}$. Otra manera de entender a los módulos p-persistentes es en términos de inclusiones $K^i \subset K^{i+p}$ con los homomorfismos inducidos $f_k^{i,p}: H_k^i \to H_k^{i+p}$ y las relaciones $im(f_k^{i,p}) \cong H_k^{i,p}$ o como los términos de una secesión espectral. (Zomorodian A. and Carlsson G., 2005)

Para conocer la información de generación y persistencia de las clases homológicas es necesario contestar las preguntas de cuándo es que estas nacen, mueren y cuánto tiempo es que viven, para ello se consideran los grupos de homología persistente que contienen clases homológicas estables en el intervalo $(i, j + 1)$, es decir que nos interesan aquellas que nacen en un tiempo no posterior a $i$ y siguen vivas en $j + 1$, ya que si permanecen para grandes valores, detectan características topológicas estables en la nube de datos X. En contraparte, para clases que permanecen únicamente para valores pequeños, estas son inestables y sus componentes topológicas no aportarán información de calidad de X, sino ruido.

**Def.2.15.** Sea K un complejo y F una filtración, la clase $\alpha \in H_p(K_i)$ :

- Nace en el tiempo $i$ si $\alpha \notin H_{i-1,i;p}(K, F)$, para las clases cero, el cero nace en $-\infty$
- Será un ancestro de la clase $\beta \in H_p(K_j)$ con $i \leq j$ si $l_{i,j}(\alpha) = \beta$, un primer ancestro es, un ancestro de $\alpha$ cuyo tiempo de nacimiento es mínimo, denotado
- Nace en un tiempo $i$ si $n(\alpha) = i$
- Muere en tiempo $j + 1$ si $l_{j,j+1}(\alpha) \in H_{i-1,j}$

**Def. 2.16.** La clase $\beta \in H_p(K_j)$ es descendiente de la clase $\alpha \in H_p(K_i)$ con $i \leq j$ si $l_{i,j}(\alpha) = \beta$, así el primer descendiente de una clase $\alpha \in H_p(K_i)$, si existe se le lama último descendiente de $\alpha$ y y a su tiempo d emuerte $m(\alpha)$, si no existe $m(\alpha) = \infty$, así el último descendiente de $\alpha$ es $l_{i,m}(\alpha)$.

**Def. 2.17.** Una colección de clases $F = \{\alpha_i, \alpha_{i+1}, \dots \alpha_j\}$ es una familia si se cumplen las siguientes condiciones:

- $l_{i,s}(\alpha_i) = \alpha_s$ para $s = i + 1, \dots j$,
- $\alpha_i$ es primer ancestro de $\alpha_j$,
- $\alpha_j$ es último ancestro descendiente de $\alpha_i$

**Def. 2.18.** Sea $\alpha \in H_p(K_i)$ una clase no cero, su persistencia se define como $pers(\alpha) = m(\alpha) - n(\alpha)$, si $m(\alpha) = \infty$, entonces $pers(\alpha) = \infty$.

**Def. 2.19**. Si $F = \{\alpha_i, \alpha_{i+1}, \dots \alpha_j\}$ es una familia, su persistencia se define como $pers(F) =$

$pers(\alpha_j)$ si $m(\alpha_j) = \infty$, entonces $pers(F) = \infty$, es claro que por ser una familia $pers(F) = pers(\alpha)$ para cualquier $\alpha \in F$.

**Def. 2.20.** Una familia $F = \{\alpha_i, \alpha_{i+1}, \dots \alpha_j\}$ es de generación $(i,j)$ y una clase $\alpha$ es de generación $(i,j)$ si pertenece a una familia de generación $(i,j)$, así una clase puede pertenecer a distintas familias pero todas estas tienen la misma generación, es decir que está bien definida.

La persistencia nos dice qué tan antigua es una clase, persistencias pequeñas no serán de gran interés para su estudio pues no representan cualidades importante en el crecimiento del complejo K, sino como se ha mencionado pueden indicar la presencia de ruido, por otro lado las persistencias grandes o que permanecen mucho tiempo y no mueren, aportarán información de gran importancia sobre K.

**Diagrama de persistencia**

Como se mencionó anteriormente, para visualizar la evolución topológica de la homología persistente en el tiempo con respecto a el parámetro $\varepsilon > 0$, se utiliza una representación gráfica conocida como diagrama de persistencia, el cual expresa la cantidad y estabilidad de los diferentes huecos k-dimensionales presentes en X, para cada nivel o generación de la homología k, que por el teorema fundamental de persistencia que más adelante se enunciará, es equivalente a la información que proporcionan los números de Betti persistentes.

Para obtener un diagrama de persistencia, es necesario seguir un algoritmo a continuación descrito, basado en los conocimientos previos. (Espinosa, M.E., 2015)

- Se selecciona $0 \leq p < n$ para construir el diagrama de codificación asociado a la homología p-dimensional

- Se elije $1 \leq i < j \leq \infty$ y se marca el plano $\mathbb{R} \times (\mathbb{R} \cup \{\infty\})$, el punto $(i, j+1)$ donde $n+1$ se marca en $\infty$ si $\mu_{i,j} > 0$

- Se le asigna multiplicidad $\mu_{i,j}$ al punto $(i, j+1)$, si este punto fue marcado y esta se marca en el diagrama escribiendo a $\mu_{i,j}$ al lado del punto $(i, j+1)$

- Se dibuja la diagonal $x = y$ además de considerar que cada punto tiene multiplicidad infinita.

La *Figura 6* ilustra la manera en que se ve un diagrama de persistencia.

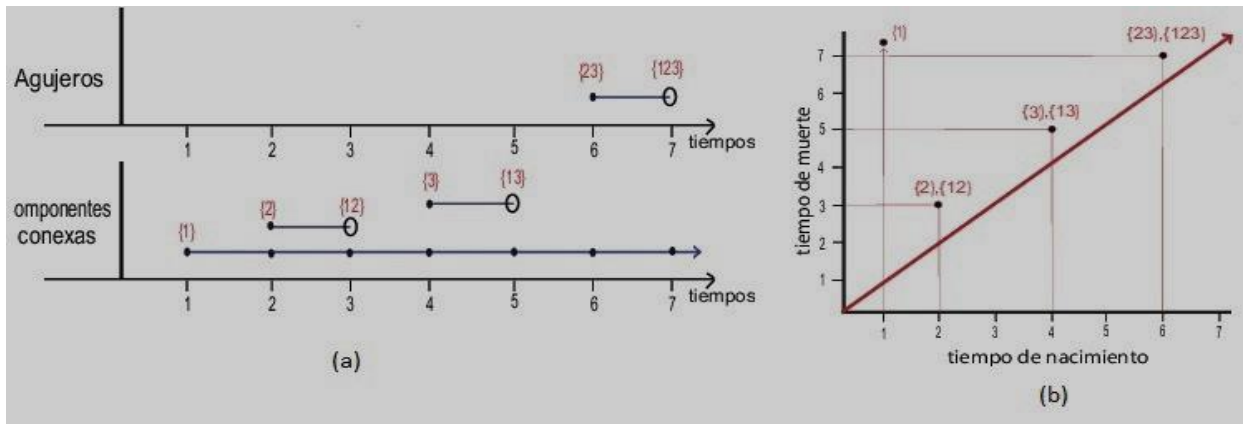

**Figura 6.** Diagramas para representar la persistencia homológica de la filtración de la figura… En (a) se representan los diagramas códigos de barra, uno para cada dimensión. En (b) se representa el llamado diagrama de persistencia. (Alonso R., Garcia E. y Lamar J., 2015, p.22)

## Conclusiones

Cada día se vuelve una tarea de mayor interés el análisis e interpretación de datos para distintas empresas, tradicionalmente el uso de minería de datos, modelos estadísticos, aproximaciones bayesianas son las más comunes, la topología algebraica aplicada a datos, es una alternativa novedosa, y que para manejo de grandes volúmenes de datos de forma visual lo hace sencillo con un rigor y confiabilidad necesaria para saber que es y será una herramienta proba, este tipo de análisis es de análisis interesante para temas de prevención de lavado de dinero, prevención de fraudes fiscales, alertas tempranas en datos masivos, vinculación de internet de las cosas con desempeño operativo.

## Referencias